\documentclass[11pt]{amsart}

\usepackage[frame,cmtip,arrow,matrix,line,graph,curve]{xy}
\usepackage{amsmath,amssymb, amscd}
\usepackage{amsfonts,amsxtra,amsthm}
\usepackage{mathrsfs}
\usepackage{eucal}
\usepackage{latexsym}

\numberwithin{equation}{section}

\newtheorem{theorem}{Theorem}[section]
\newtheorem{proposition}[theorem]{Proposition}
\newtheorem{corollary}[theorem]{Corollary}
\newtheorem{lemma}[theorem]{Lemma}

\newtheorem{problem}[theorem]{Problem}
\newtheorem{remit}[theorem]{Remark}
\newtheorem{definit}[theorem]{Definition}

\newenvironment{definition}{\begin{definit}\rm}{\end{definit}}

\newenvironment{remark}{\begin{remit}\rm}{\end{remit}}

\def\lab{\label}

\newcommand{\pp}{\mathbb{P}}

\newcommand{\cc}{\mathbb{C}}

\newcommand{\zz}{\mathbb{Z}}

\newcommand{\Hom}{\mathrm{Hom}}

\newcommand{\git}{/\!\!/}

\newcommand{\cA}{\mathcal{A} }
\newcommand{\cB}{\mathcal{B} }

\newcommand{\cE}{\mathcal{E} }
\newcommand{\cH}{\mathcal{H} }
\newcommand{\cI}{\mathcal{I} }
\newcommand{\cK}{\mathcal{K} }
\newcommand{\cM}{\mathcal{M} }
\newcommand{\cN}{\mathcal{N} }
\newcommand{\cL}{\mathcal{L} }
\newcommand{\cO}{\mathcal{O} }

\newcommand{\cU}{\mathcal{U} }

\newcommand{\cW}{\mathcal{W} }

\newcommand{\tX}{\tilde{X} }

\newcommand{\tQ}{\widetilde{Q} }
\newcommand{\hQ}{\widehat{Q} }
\newcommand{\tcM}{\widetilde{\cM} }
\newcommand{\hcM}{\widehat{\cM} }

\newcommand{\s}{\sigma }

\newcommand{\bM}{\mathbf{M} }

\newcommand{\bH}{\mathbf{H} }
\newcommand{\bC}{\mathbf{C} }

\newcommand{\obM}{\overline{\bM} }
\newcommand{\obMz}{\overline{\bM} _{0,0} }

\begin{document}

\title[Stable maps to moduli space]
{Hecke correspondence, stable maps and the Kirwan
desingularization}
\date{November 30, 2005}

\author{Young-Hoon Kiem}
\address{Department of Mathematics and Research Institute
of Mathematics, Seoul National University, Seoul 151-747, Korea}
\address{Department of Mathematics, Stanford University, Stanford,
USA} \email{kiem@math.snu.ac.kr}
\thanks{Partially supported by KRF grant 2005-070-C00005}

\begin{abstract}
We prove that the moduli space $\obMz(\cN,2)$ of stable maps of
degree 2 to the moduli space $\cN$ of rank 2 stable bundles of
fixed determinant over a smooth projective curve of genus $g$ has
two irreducible components which intersect transversely. One of
them is Kirwan's partial desingularization $\widetilde{\cM}_X$ of
the moduli space $\cM_X$ of rank 2 semistable bundles with
determinant isomorphic to $\cO_X(y-x)$ for some $y\in X$. The
other component is the partial desingularization of $\pp \Hom
(sl(2)^\vee, \cW)\git PGL(2)$ for a vector bundle
$\cW=R^1{\pi}_*\cL^{-2}(-x)$ of rank $g$ over the Jacobian of $X$.
We also show that the Hilbert scheme $\bH$, the Chow scheme $\bC$
of conics in $\cN$ and $\obMz(\cN,2)$ are related by explicit
contractions.
\end{abstract}
\maketitle


\section{Introduction}

Let $X$ be a smooth projective curve of genus $g\ge 3$ over the
complex number field. We fix $x\in X$ throughout this paper. The
moduli space $\cN$ of rank 2 stable bundles over $X$ with
determinant $\cO_X(-x)$ is a smooth projective Fano variety whose
Picard group is generated by a very ample line bundle $\Theta$.
(See \cite{BV}.) The problem that we address in this paper is the
following.

\begin{problem}\lab{prob0.1} Describe explicitly the moduli space
$\obMz(\cN,d)$ of stable maps to $\cN$ of genus 0 and degree
$d$.\end{problem}

Obviously the degree of a stable map $f:C\to \cN$ is defined as
the degree of $f^*\Theta$. The expected dimension of
$\obMz(\cN,d)$ is $3g-6+2d$ by Riemann-Roch since
$c_1(T\cN)=2\Theta$. When $d=1$, V. Munoz \cite{Munoz} proved that
$\obMz(\cN,1)$ is isomorphic to a Grassmannian bundle $Gr(2,\cW)$
over $J=Pic^0(X)$. (See also \cite{Sun}.) Here
$\cW=R^1{\pi_J}_*\cL^{-2}(-x)$ is a vector bundle of rank $g$ over
$J$ where $\cL$ is a Poincar\'e bundle over $J\times X$ and
$\pi_J:J\times X\to J$ is the projection.

The purpose of this paper is to study the $d=2$ case. This is
particularly interesting because of Hecke correspondence which has
been one of the most powerful tools in the study of moduli spaces
of bundles. Let $\cM_X$ be the moduli space of rank 2 semistable
bundles $E$ with $\det E\cong \cO_X(y-x)$ for some $y\in X$. This
is a $(3g-2)$-dimensional singular projective normal variety which
fits into a Cartesian diagram
\[\xymatrix{
\cM_X\ar@{^(->}[r]\ar[d] & \cM(2,0)\ar[d]^{\det}\\
X\ar@{^(->}[r] & J }\] where $M(2,0)$ is the moduli space of
semistable bundles of rank 2 and degree 0 over $X$ and the bottom
map is $y\mapsto \cO_X(y-x)$. Hecke correspondence refers to a
diagram
\begin{equation}\lab{eq0.1}
\xymatrix{ &\cH\ar[dl]_{q_1}\ar[dr]^{q_0}\\
\cN\times X && \cM_X}\end{equation} where $q_1$ is the projectivization
of a universal bundle over $\cN\times X$ and $q_0$ is a
$\pp^1$-bundle over the stable part $\cM_X^s$. Hence for any
$\theta\in \cM_X^s$, we have a stable map
 \[ \xymatrix{ q_0^{-1}(\theta)\ar[r]^{q_1} &\cN\times
 X\ar[r]^{\pi_\cN} &\cN }\]
which is of degree 2 (\cite{NR}). By \eqref{eq0.1}, $\cM_X^s$
parameterizes stable maps to $\cN$ and hence we have a morphism
$$\cM_X^s\to \obMz(\cN,2)$$
which is injective and birational onto an irreducible component.
We prove that this component is precisely Kirwan's partial
desingularization $\tcM_X$ of $\cM_X$.

Kirwan's partial desingularization is a systematic way to
partially resolve the singularities of GIT quotients. The partial
desingularization $\tcM_X$ of $\cM_X$ (see \cite{k5}) is the consequence of
two blowups first along the deepest stratum
\begin{equation}\lab{eq0.2}
\tX:=\{\xi\in J\,|\, \xi^2\cong \cO_X(y-x)\text{ for some }y\in
X\}
\end{equation}
which is the locus of bundles $\xi\oplus \xi$, $\xi\in \tX$ and
then along the proper transform of the middle stratum
\begin{equation}\lab{eq0.3}
\cK_J:=\{(\xi_1,\xi_2)\in J\times J\,|\,\xi_1\xi_2\cong
\cO_X(y-x)\text{ for some }y\in X\}/\zz_2\end{equation} which is
the locus of bundles $\xi_1\oplus \xi_2$, $(\xi_1,\xi_2)\in \cK_J$
where $\zz_2$ interchanges $\xi_1$ and $\xi_2$. The upshot is an
orbifold $\tcM_X$ which is singular along a nonsingular subvariety
$\Delta_{\tX}\cong \pp (S^2\cA^\vee_{\tX})$ where $\cA_{\tX}\to
Gr(2,\cW_0)$ is the universal rank 2 bundle on the relative
Grassmannian of the bundle $\cW_0=R^1{\pi}_*(\cL_0^{-2}(-x))$,
$\cL_0$ being the Poincar\'e line bundle over $\tX\times X$ and
$\pi:\tX\times X\to\tX$ being the projection. By blowing up
$\tcM_X$ along $\Delta_{\tX}$ we obtain a (full) desingularization
$\hcM_X$ of $\cM_X$. The exceptional divisor of the last blowup is
\begin{equation}\lab{eq0.4}
\widehat{\Delta}_{\tX}=\pp
(S^2\cA_{\tX}^\vee)\times_{Gr(2,\cW_0)}\pp(\cW_0/\cA_{\tX}\oplus
\eta)
\end{equation} for some line bundle $\eta$ over $Gr(2,\cW_0)$
and the normal bundle is $\cO(-1,-2)$ on the fibers $\pp^2\times
\pp^{g-2}$ over $Gr(2,\cW_0)$. The blowup $\hcM_X\to\tcM_X$ is the
contraction of $\pp (\cW_0/\cA_{\tX}\oplus \eta)$.

Although the partial desingularization $\tcM_X$ has been quite
useful especially for cohomological computations, its moduli
theoretic meaning has been unknown. We prove the following in this
paper.
\begin{theorem}\lab{thm0.1}
(1) The partial desingularization $\tcM_X$ is isomorphic to the
irreducible component of $\obMz(\cN,2)$ containing the image of
$\cM_X^s$.

(2) We have \begin{equation}\lab{eq0.5} \obMz(\cN,2)\cong
\tcM_X\cup \tQ_J\end{equation} where $\tQ_J=\obMz(\pp\cW/J,2)$ is
the partial desingularization of the GIT quotient $Q_J=\pp \Hom
(sl(2)^\vee,\cW)\git PGL(2)$.
 The two irreducible components
intersect transversely along the partial desingularization
$\tQ_{\tX}=\obMz(\pp\cW_0/\tX,2)$ of the GIT quotient
$Q_{\tX}=\pp \Hom (sl(2)^\vee,\cW_0)\git PGL(2)$.
\end{theorem}

In \cite{JKKW1}, a strategy for computing the intersection pairing
of the partial desingularization of a GIT quotient is provided and
the intersection pairing of the partial desingularization of the
moduli space $\cM(2,\cO_X)$ was studied in \cite{JKKW2}. By
Theorem \ref{thm0.1}, the results of \cite{JKKW2} should give us
the second order deformation of the quantum cohomology ring of
$\cN$. Since the quantum cohomology ring is in principle
understood \cite{Munoz}, we didn't pursue the computation.

A problem closely related to Problem \ref{prob0.1} is the
following. \begin{problem} Describe explicitly the Hilbert scheme
and the Chow scheme of curves of degree $d$, i.e. subschemes of
$\cN$ with Hilbert polynomial $dm+1$ with respect to the ample
generator $\Theta$.\end{problem} When $d=1$, they all coincide
with $\obMz(\cN,1)$. We prove the following for conics (degree 2
curves).

\begin{theorem}\lab{thm0.2}
The Hilbert scheme $\bH$ of conics, the Chow scheme $\bC$ of
conics and $\bM:=\obMz(\cN,2)$ are related by contractions
\begin{equation}\lab{eq0.6}
\xymatrix{ & \widehat{\bM}\ar[dl]
\ar[dr]^{\pp(S^2\cA_{J}^\vee)}\\
\bM\ar[dr]_{\pp(S^2\cA_{J}^\vee)}
&&\bH\ar[dl]\\
&\bC }\end{equation} where $\widehat{\bM}$ is the blowup of the
disjoint union of the two components of $\bM$ along the singular
loci. The two irreducible components of $\bH$ are nonsingular.
\end{theorem}

This paper is organized as follows. In \S2 and \S3, we classify
rational curves in $\cN$. In \S4, we study the moduli space of
stable maps to a projective space. In \S5, we recall the Hecke
correspondence and the partial desingularization $\tcM_X$. We
prove Theorem \ref{thm0.1} in \S6 and Theorem \ref{thm0.2} in \S7.

I am most grateful to Professor Michael Thaddeus for sharing his
remarkable insight that the partial desingularization $\tcM_X$ must be the
irreducible component of $\obMz(\cN,2)$ containing Hecke curves. I
am also grateful to Professor Ravi Vakil for answering a few
questions. This paper was written while the author was visiting
the mathematics department of Stanford University. It is my
pleasure to express gratitude to Professor Jun Li and the
mathematics department of Stanford for generous support and
hospitality during my stay.

After finishing this paper, Professor Jun-Muk Hwang kindly
informed me of \cite{Cas} where A.-M. Castravet studied rational
curves in $\cN$. There seems to be some overlap with her paper and
\S3 of this paper. I apologize for any overlap with \cite{Cas}
although our motivations and main results are totally independent.


\section{Examples of stable maps to $\cN$}
In this section we consider various examples of rational curves in
the moduli space $\cN$ of stable rank 2 bundles with determinant
$\cO_X(-x)$.
\subsection{Hecke curves}
Let $E\in \cM_X^s$ be a stable rank 2 bundle with $\det
E\cong \cO_X(y-x)$ for some $y\in X$. For $\nu\in \pp
E|_y^\vee\cong \pp^1$, let $E^\nu$ denote the kernel of the
composition
\[ \xymatrix{ E\ar[r] & E|_y\ar[r]^\nu & \cc}\]
of $\nu$ and the restriction to $y$. Then $E^\nu$ is a stable
bundle with $\det E^\nu\cong \cO(-x)$. On $\pp E|_y\times
X=\pp^1\times X$, we have a canonical map
\begin{equation}\lab{eq3.1}
\pi_X^*E\to \pi_X^*E|_y\to \cO_{\pp E|_y\times\{y\}}(1)
\end{equation}
where $\pi_X$ denotes the projection onto $X$. The kernel $\cE$ of
\eqref{eq3.1} is a family of rank 2 bundles on $X$ with
determinant $\cO(-x)$, parameterized by $\pp^1$. Therefore we
obtain a morphism $f_{E,y}:\pp^1\to \cN$ since $\cN$ is the moduli
space of such bundles.
\begin{definition} For $E\in \cM_X^s$ with  $\det E\cong \cO_X(y-x)$,
the induced morphism $f_{E,y}:\pp^1\to \cN$ is called a
\emph{Hecke curve}.\end{definition} Conversely, suppose the
pull-back $\cE$ of a universal bundle on $\cN\times X$ to
$\pp^1\times X$ by $f\times id_X$ for a morphism $f:\pp^1\to \cN$
fits into an exact sequence
\begin{equation}\lab{eq3.2}
0\to \cE\to \pi_X^*E\to \cO_{\pp^1\times\{y\}}(1)\to 0.
\end{equation} Then $f$ is a Hecke curve as it is given by
$E|_y\otimes \cO_{\pp^1\times\{y\}}\to \cO_{\pp^1\times\{y\}}(1)$
which is unique modulo the action of $PGL(2)$.

Upon restricting \eqref{eq3.2} to $\pp^1\times \{y\}$ we get a map
\[ \cE|_{\pp^1\times\{y\}}\to \cO_{\pp^1\times\{y\}}(-1)\to 0\]
and thus an exact sequence
\begin{equation}\lab{eq3.3}
0\to \pi_X^*E(-y)\to \cE\to \cO_{\pp^1\times\{y\}}(-1)\to
0.\end{equation} So we obtain the following with $F=E(-y)$.
\begin{lemma}
Let $f:\pp^1\to \cN$ be a morphism such that the pull-back $\cE$
of the universal bundle on $\cN\times X$ by $f\times id_X$ fits
into an exact sequence \begin{equation}\lab{eq3.4} 0\to
\pi_X^*F\to \cE\to \cO_{\pp^1\times\{y\}}(-1)\to 0 \end{equation}
for a rank 2 bundle $F$ and $y\in X$. Then $f$ is a Hecke curve.
\end{lemma}

We recall the following results of Narasimhan and Ramanan.
\begin{theorem} \lab{thm3.100} \cite[\S5]{NR}
\begin{enumerate}
\item For any stable bundle $E$ with $\det E\cong \cO_X(y-x)$, $f_{E,y}$ is an
embedding and the image is a smooth rational curve in $\cN$ of
degree 2. \item The normal bundle of a Hecke curve is generated by
global sections. \item For distinct $E\in \cM_X^s$, the
corresponding Hecke curves are distinct.\end{enumerate}
\end{theorem}
\begin{proof}
We only provide a proof of (3). Let $f:\pp^1\to \cN$ be a Hecke
curve and let $\cE\to \pp^1\times X$ be the pull-back of a
universal bundle over $\cN\times X$ such that $\cE|_{\pp^1\times
z}\cong \cO\oplus \cO$ for general $z \in X$. Then
$E={\pi_X}_*\cE\otimes \cO_X(y)\in \cM_X^s$ is the unique rank 2
bundle which gives us the Hecke curve $f$ by Theorem \ref{thm4.1}
and \eqref{eq3.3}.
\end{proof}

\begin{remark} \lab{rmk2.4n} Item (2) implies that $H^1(\pp^1, f^*T\cN)=0$ and
hence $\obMz(\cN, 2)$ is nonsingular at the points of Hecke
curves.
\end{remark}

With the Hecke curves, we can generate many more rational curves
of higher degrees. First, we can compose a degree $l$ map
$\pp^1\to \pp^1$ with a Hecke curve $f_{E,y}:\pp^1\to \cN$ and get
a degree $2l$ map $f:\pp^1\to \cN$. In this case the pull-back
$\cE$ of the universal bundle fits into an exact sequence
\begin{equation}\lab{eq3.5}
0\to \pi_X^*F\to \cE\to \cO_{\pp^1\times\{y\}}(-l)\to 0
\end{equation} where $F=E(-y)$. Conversely, if a morphism
$f:\pp^1\to \cN$ gives us \eqref{eq3.5}, it should be the
composition of a Hecke curve with a degree $l$ covering $\pp^1\to
\pp^1$.

More generally, let $l_1,\cdots, l_r$ be a sequence of positive
integers and $y_1,\cdots, y_r$ be distinct points in $X$. Let $E$
be a rank 2 bundle on $X$ with $\det E\cong \cO(y_1+\cdots +y_r
-x)$ such that \begin{equation}\lab{eq3.5a}\deg L < \frac12 (\deg
E -r)\text{ for any line subbundle }L.\end{equation} Over $\pp
E|_{y_1}^\vee\times \cdots \times \pp E|_{y_r}^\vee\times X$, we
have a canonical map
\[ \pi_X^*E\to \pi_X^*E|_{y_1}\oplus \cdots \oplus \pi_X^*E|_{y_r}\to \cO_{\pp
E|_{y_1}^\vee\times\{y_1\}}(1)\oplus \cdots \oplus \cO_{\pp
E|_{y_r}^\vee\times\{y_r\}}(1) \] whose kernel is a family of
stable bundles with determinant $\cO_X(-x)$. Thus we have a
morphism
\[ \pp E|_{y_1}^\vee\times \cdots \times\pp
E|_{y_r}^\vee\to\cN.
\]
Now degree $l_i$ maps $\pp^1\to \pp E|_{y_i}^\vee$ give us a map
\[ f:\pp^1\to \pp E|_{y_1}^\vee\times \cdots \times\pp
E|_{y_r}^\vee\to\cN
\]
whose degree is $2(l_1+\cdots +l_r)$. Let us call this a
\emph{generalized Hecke curve}. The pull-back $\cE$ of a
universal bundle fits into an exact sequence
\begin{equation}\lab{eq3.6}
0\to \pi_XF\to \cE\to \oplus_i \cO_{\pp
E|_{y_i}^\vee\times\{y_i\}}(-l_i)\to 0
\end{equation} where $F=E(-y_1-\cdots -y_r)$. Conversely, if
$f:\pp^1\to \cN$ induces $\cE$ given by \eqref{eq3.6}, then $f$ is
a generalized Hecke curve.
\begin{remark} \eqref{eq3.5a} requires that every elementary modification by
an element of $\pp E|_{y_1}^\vee\times \cdots \times\pp
E|_{y_r}^\vee$ be stable. But to get a generalized Hecke curve we
only have to require that the elementary modifications by elements
in the image of $f$ be stable. \end{remark}

\subsection{Extensions of line bundles}

A different way of producing rational curves in $\cN$ is to
consider extension bundles. For $\xi\in Pic^k(X)$, let $\pp
Ext^1(\xi, \xi^{-1}(-x))^s=\pp H^1(\xi^{-2}(-x))^s$ be the locus
of extension bundles
\[ 0\to \xi^{-1}(-x)\to E\to \xi\to 0\]
which are stable. Then we have an obvious morphism \[ \theta_\xi:
\pp Ext^1(\xi, \xi^{-1}(-x))^s\to \cN.\] Suppose $f':\pp^1\to \pp
Ext^1(\xi,\xi^{-1}(-x))$ is a morphism of degree $a$ whose image
lies in $\pp Ext^1(\xi, \xi^{-1}(-x))^s$. Then the composition of
$f'$ with $\theta_\xi$ gives us a stable map $f:\pp^1\to \cN$.

When $k=0$ it is easy to check that $\pp Ext^1(\xi,
\xi^{-1}(-x))^s=\pp Ext^1(\xi, \xi^{-1}(-x))$. We recall the
following.
\begin{proposition}\lab{prop2.134}
\cite[6.18]{NR} $\theta_\xi$ is an embedding of degree 1 for
any $\xi\in Pic^0(X)$, i.e. the pull-back of the ample generator
of $Pic(\cN)$ is $\cO(1)$. \end{proposition}

In general, we have the following.
\begin{lemma}\lab{lem2.117} For $\xi\in Pic^k(X)$, the degree of $\theta_\xi$ is
$2k+1$.\end{lemma}
\begin{proof}
The universal extension bundle $\cU$ over $\pp Ext^1(\xi,
\xi^{-1}(-x))^s\times X$ fits into an exact sequence
\[ 0\to\pi_X^*\xi^{-1}(-x)\otimes \cO_{\pp}(1)\to \cU\to
\pi_X^*\xi\to 0\] By \cite[6.7]{NR}, the tangent bundle $T\pp$ of
$\pp Ext^1(\xi, \xi^{-1}(-x))^s$ fits into an exact sequence
\[ 0\to T\pp\to R^1{\pi_\pp}_*(\cU\otimes \pi_X^*\xi^{-1})\to
R^1{\pi_\pp}_*\cO=H^1(X,\cO)\otimes \cO_{\pp}\to 0\] where
$\pi_\pp$ denotes the projection onto $\pp Ext^1(\xi,
\xi^{-1}(-x))^s$. Since $$\dim \pp Ext^1(\xi,
\xi^{-1}(-x))^s=2k+g-1$$ by Riemann-Roch,
\[ c_1(R^1{\pi_\pp}_*(\cU\otimes
\pi_X^*\xi^{-1}))=c_1(T\pp)=\cO(2k+g).\] Also we have an exact
sequence
\[ 0\to \cU\otimes \pi_X^*\xi^{-1}\to End_0\cU\to
\pi_X^*\xi^2(x)\otimes \cO_\pp (-1)\to 0\] which gives rise to an
exact sequence
\[ 0\to H^0(\xi^2(x))\otimes \cO_\pp(-1)\to R^1{\pi_\pp}
_*(\cU\otimes \pi_X^*\xi^{-1})\to R^1{\pi_\pp} _*End_0\cU\to
H^1(\xi^2(x))\otimes \cO_\pp(-1)\to 0.\] By Riemann-Roch,
\[ \det \left( H^0(\xi^2(x))\otimes \cO_\pp(-1)
\right)^{-1}\otimes \det \left( H^1(\xi^2(x))\otimes
\cO_\pp(-1)\right) \cong \cO_\pp(2k-g+2)\] Since
$\theta_\xi^*T\cN=R^1{\pi_\pp}_*End_0\cU$,
$\deg(f^*T\cN)=(2k-g+2)+(2k+g)=4k+2.$ Because $c_1(T\cN)$ is 2
times the ample generator of $Pic(\cN)$, the degree of $f$ is
$2k+1$.
\end{proof}

As a consequence, a degree $a$ map $f':\pp^1\to \pp
Ext^1(\xi,\xi^{-1}(-x))^s$ for some $\xi\in Pic^k(X)$ composed
with $\theta_\xi$ is a map $f:\pp^1\to \cN$ of degree $a(2k+1)$.
The pull-back of a universal bundle on $\cN\times X$ by $f\times
id_X$ embeds into an exact sequence \begin{equation}\lab{eq3.40}
0\to \pi_X^*\xi^{-1}(-x)\otimes \cO_\pp(a)\to \cE\to \pi_X^*\xi\to
0.\end{equation} Conversely, it is obvious that any such $\cE$ is
given by a morphism \[ f:\pp^1\to \pp Ext^1(\xi,\xi^{-1}(-x))^s\to
\cN. \]

More generally, if we have a \emph{rational} map
$f':\pp^1\dashrightarrow \pp Ext^1(\xi,\xi^{-1}(-x))^s$ of degree
$a$, the composition of $f'$ with $\theta_\xi$ uniquely extends to
a morphism $f:\pp^1\to \cN$. The pull-back $\cE$ of the universal
bundle on $\cN\times X$ now embeds into an exact sequence
\begin{equation}\lab{eq3.41}
0\to \pi_X^*\xi^{-1}(-x)\otimes \cO_\pp(a)\to \cE\to
\pi_X^*\xi\otimes \cI_Z\to 0\end{equation} for a zero-dimensional
subscheme $Z$ of $\pp^1\times X$. (See Theorem \ref{thm4.1} (2).)


\section{Classification of rational curves in $\cN$}

In this section we show that the examples discussed in \S3 are all
possible nonconstant maps $f:\pp^1\to \cN$.

Suppose $f:\pp^1\to\cN$ is a nonconstant morphism. The pull-back
of a universal bundle over $\cN\times X$ by $f\times id_X$ is a
rank 2 bundle
\[  \cE\to \pp^1\times X=:S\]
As before, let $\pi_X,\pi_P$ denote the projections of $S$ to $X$
and $\pp^1$ respectively. The starting point of our classification
is the following lemma due to X. Sun.
\begin{lemma}\lab{lem4.00}
\cite[\S2]{Sun}
Let $f:\pp^1\to \cN$ be a morphism of degree $\beta\in \zz_{>0}$
and $\cE$ be the induced vector bundle. Then
$$\beta=2c_2(\cE)-\frac12 c_1(\cE)^2.$$
\end{lemma}
\begin{proof}
It is well-known that $f^*T\cN=R^1{\pi_P}_*End_0\cE$ and
${\pi_P}_*End_0\cE=0$ by stability where $End_0\cE$ is the
traceless part of the endomorphism bundle of $\cE$. By Leray's
spectral sequence and Riemann-Roch, we have
\[ \chi(S,End_0\cE)=-\chi(\pp^1,f^*T\cN)=-(3g-3+2\beta)\]
since $\deg c_1(f^*T\cN)=\deg f^*c_1(T\cN)=2\beta$. On the other
hand, by Riemann-Roch for $S$ we have
\[\chi(S,End_0\cE)=\int_S ch(End_0\cE)\cdot td_X\cdot
td_{\pp^1}=-3(g-1)-c_2(End_0\cE).\] Therefore $\beta=\frac12
c_2(End_0\cE)=2c_2(\cE)-\frac12 c_1(\cE)^2.$ \end{proof}

By tensoring with $\pi_P^*\cO(d)$ for suitable $d\in \zz$ we may
assume that $\cE|_{\pp^1\times \{z\}}$ for general $z$ is
$\cO_{\pp^1}(a)\oplus \cO_{\pp^1}$ for some $a\ge 0$. This doesn't
change the morphism $f$ since tensoring by $\pi_P^*\cO(d)$ gives
us an equivalent family of stable bundles. We recall the following
result of Brosius.

\begin{theorem} \lab{thm4.1} \cite[Theorem 1]{Bros}
Let $\cE$ be as above. Then $\cE$ belongs to an exact sequence
\[
0\to \cE'\to \cE\to \cE''\to 0
\]
such that the following hold. \begin{enumerate} \item If $a=0$,
there exist a rank 2 bundle $F$ on $X$ and $Z\in Hilb^l(S)$ for
some $l\ge 0$ such that $\cE'=\pi_X^*F$ and $\cE''=\cI_{Z\subset
Y}$ is the ideal sheaf of $Z$ in $Y=\pi_X^{-1}\pi_X(Z)$. \item If
$a>0$, there exist a line bundle $\xi\in Pic^k(X)$ for some $k\ge
0$ and $Z\in Hilb^l(S)$ for some $l\ge 0$ such that
$\cE'=\pi_X^*\xi^{-1}(-x)\otimes \pi_P^*\cO_{\pp^1}(a)$ and
$\cE''=\cI_Z\otimes \pi_X^*\xi$.
\end{enumerate}
Furthermore, $\cE'$ and $\cE''$ are unique up to isomorphism.
\end{theorem}

\begin{definition} We say $f$ is of \emph{Hecke type} when (1) holds and
$f$ is of \emph{extension type} when (2) holds.
\end{definition}

\begin{remark}
(1) Suppose $f$ is of Hecke type. If $l=1$, $f$ is a Hecke curve
by \eqref{eq3.4} because $\cI_{Z\subset Y}\cong \cO_{\pp^1\times
\{y\}}(-1)$ for $y=\pi_X(Z)$. When $l>1$, $f$ is a generalized
Hecke curve by \eqref{eq3.6}.

(2) Any extension
$$0\to \pi_X^*\xi^{-1}(-x)\otimes \pi_P^*\cO_{\pp^1}(a)\to \cE
\to \cI_Z\otimes \pi_X^*\xi\to 0$$ is locally free since the
Cayley-Bacharach property is satisfied (\cite[Theorem
5.1.1]{HL97}) from the vanishing \small
\[ H^0(\pi_X^*K_X\otimes \pi_X^*\xi^2(x)\otimes
\pi_P^*\cO(-a)\otimes \pi_P^*K_{\pp^1})=H^0(X,K_X\otimes \xi^2(x))\otimes
H^0(\pp^1,\cO(-a-2))=0.\]\normalsize

(3) If $f$ is of extension type and $l=0$, then $f$ is a line in
$\pp H^1(\xi^{-2}(-x))^s$. When $l>0$, $f$ is given by a rational
map $\pp^1\dashrightarrow \pp H^1(\xi^{-2}(-x))^s$ of degree $a$.
\end{remark}

By direct computation with Lemma \ref{lem4.00},
we obtain the following useful lemma.
\begin{lemma}\lab{keylem} With notation as in Theorem \ref{thm4.1},
$\beta=2l$ for the Hecke
type maps and $\beta=2l+a(2\deg \xi+1)$ for the extension type
maps.\end{lemma}

We can now classify rational curves of degree $\le 2$ in $\cN$.
\begin{proposition}\lab{prop4.1} Let $f:\pp^1\to \cN$ is of degree $\le 2$.
Then $f$ is one of the following.
\begin{enumerate}
\item $f$ is a Hecke curve for some $y\in X$ and a stable rank 2
bundle $E$ with determinant $\cO_X(y-x)$ \item $f$ is a conic in
$\pp H^1(\xi^{-2}(-x))$ for some $\xi\in Pic^0(X)$ \item $f$ is a
line in $\pp H^1(\xi^{-2}(-x))$ for some $\xi\in
Pic^0(X)$.\end{enumerate} The degree of $f$ is $2$ for (1) and (2)
while $1$ for (3). \end{proposition}
\begin{proof}
By Lemma \ref{keylem}, $\beta=1$ is possible only when $f$ is of
extension type with $l=0$, $a=1$, $\deg \xi=0$. This holds exactly
when (3) above holds. Suppose $\beta=2$ now. If $f$ is of Hecke
type, then $l=1$ which means $f$ is a Hecke curve. If $f$ is of
extension type, then $l=0$, $a=2$ and $\deg\xi=0$ which means that
$f$ is a conic in $\pp H^1(\xi^{-2}(-x))$. \end{proof}

As a first consequence, we obtain the moduli space $\obMz(\cN,1)$
of stable maps of degree 1.
\begin{theorem} \cite{Munoz, Sun} \lab{thm3.18}
Let $\cL\to Pic^0(X)\times X$ be a Poincar\'e bundle and let
$\pi:Pic^0(X)\times X\to Pic^0(X)$ be the projection. Let $$Gr(2,
\cW)\to Pic^0(X)=J$$ be the Grassmannian bundle of 2 dimensional
subspaces of the bundle
$$\cW=R^1\pi_*\cL^{-2}(-x)=R^1\pi_*Hom(\cL, \cL^{-1}(-x))$$ of rank $g$
over $J$. Then we have an isomorphism \[\obMz(\cN,1)\cong Gr(2,
\cW).
\]
\end{theorem}
\begin{proof} Let $f\in \obMz(\cN,1)$.
By Proposition \ref{prop4.1}, $f$ factors through a degree 1 map
$\pp^1\to \pp H^1(\xi^{-2}(-x))$. Obviously the Grassmannian
bundle $Gr(2, \cW)$ parameterizes all such maps and thus we get a
natural map
\[ \gamma: Gr(2, \cW)\to \obMz(\cN,1)\]
which is bijective by Proposition \ref{prop2.134}
and Proposition \ref{lem1.5}.
For a degree 1 extension type map $f$, we will see in
\S6 (Proposition \ref{prop5.10}) that
$$f^*T\cN\cong \cO(2)\oplus \cO(1)^{k}\oplus \cO(-1)^{k}\oplus
\cO^l$$ for $k,l$ such that $2k+l+1=3g-3$. Therefore
$H^1(\pp^1,f^*T\cN)=0$ and thus $\obMz(\cN,1)$ is normal. To
conclude that $\gamma$ is an isomorphism by Zariski's main
theorem, it suffices to show that $\gamma$ is birational. Locally
in usual topology on the smooth part, $\gamma$ is given by
holomorphic maps $\gamma:\cc^n\to\cc^n$. If the Jacobian
determinant is nonzero at some point, then $\gamma$ is birational.
It is an elementary consequence of inverse function theorem that
if the Jacobian determinant is everywhere zero the map cannot be
one to one. So we are done.
\end{proof}

Note that $Gr(2,\cW)$ is the moduli space of relative stable maps
$\obMz(\pp \cW/J,1)$.

The open dense subset in $\obMz(\cN,2)$ of irreducible curves
consists of stable maps given in (1) and (2) of Proposition
\ref{prop4.1}. In \S6, we will prove that $\obMz(\cN,2)$ consists
of two irreducible components: Hecke curves will give us an
irreducible component of dimension $3g-2$ (isomorphic to the
partial desingularization of $\cM_X$) and the extension type maps
will give us an irreducible component of dimension $4g-4$
(isomorphic to $\tQ_J=\obMz(\pp \cW/J,2)$ in \S4).

When the degree of $f$ is 3, Lemma \ref{keylem} says that $f$ has
to be of extension type and there are three possibilities:
\begin{enumerate}
\item $l=0, a=3, \deg \xi =0$ \item $l=0, a=1, \deg \xi=1$ \item
$l=1, a=1, \deg \xi=0$.\end{enumerate} The first case is a cubic
curve in $\pp Ext^1(\xi, \xi^{-1}(-x))=\pp H^1(\xi^{-2}(-x))$ for
$\xi\in Pic^0(X)$ and the second case is a line in $\pp
H^1(\xi^{-2}(-x))^s$ for $\xi\in Pic^1(X)$. We claim the third
case is impossible.

\begin{lemma}\lab{lem4.2} Suppose $f$ is of extension type and $\deg
\xi=0$. Then $l=0$.\end{lemma} \begin{proof} Suppose $l>0$ and let
$(p,y)\in Supp Z\subset \pp^1\times X$. When $\deg \xi=0$, it is
easy to see that the restriction of $\cE$ to $\{p\}\times X$ is
unstable. Hence this is impossible.
\end{proof}

The following classification of maps $f:\pp^1\to \cN$ of degree 3
and 4 is direct from Lemmas \ref{keylem} and \ref{lem4.2}.
\begin{proposition} (1) A degree 3 map $f:\pp^1\to \cN$ comes from
either a cubic curve in $\pp H^1(\xi^{-2}(-x))$ for some $\xi\in
Pic^0(X)$ or a line in $\pp H^1(\xi^{-2}(-x))^s$ for some $\xi\in
Pic^1(X)$.

(2) A degree 4 map $f:\pp^1\to \cN$ comes from either a quartic
curve in $\pp H^1(\xi^{-2}(-x))$ for some $\xi\in Pic^0(X)$ or a
double cover of a Hecke curve. \end{proposition}


\section{Stable maps to projective spaces}

In this section we describe the moduli space
$\obM_{0,0}(\pp^{g-1},2)$ of stable maps of degree 2 in
$\pp^{g-1}$ as the partial desingularization of a GIT quotient.

Let $V$ be a vector space of dimension $g$ and let $f:\pp^1\to \pp
V=\pp^{g-1}$ be a degree 2 embedding. Then we have a surjection
\begin{equation}\lab{eq2.1}
V^\vee=H^0(\pp V,\cO(1))\twoheadrightarrow
H^0(\pp^1,\cO(2))=H^0(T\pp^1)\cong sl(2)\end{equation} because
$H^0(T\pp^1)=\mathrm{Lie} Aut(\pp^1)=sl(2)$. Conversely, any such
surjection gives us an embedding of degree 2. Hence the moduli
space of degree 2 embeddings of $\pp^1$ in $\pp V$ is the quotient
\[
\pp \mathrm{Hom}_3(sl(2)^\vee, V)/PGL(2)
\]
where the subscript 3 denotes the open locus of rank 3
homomorphisms. The compactification $\obM_{0,0}(\pp V, 2)$ turns
out to be Kirwan's partial desingularization of the GIT quotient
(\cite{k2})
\[
Q:=\pp \mathrm{Hom}(sl(2)^\vee, V)\git PGL(2).
\]

\begin{theorem}\lab{thm1.1}
Let $P_1=\pp \mathrm{Hom}_1(sl(2)^\vee, V)$ denote the locus of
rank 1 homomorphisms in $P:=\pp \Hom (sl(2)^\vee, V)$. Let
$\widetilde{P}$ be the blowup of $P$ along $P_1$ and let
$\widetilde{P}^s$ be the set of stable points in $\widetilde{P}$
with respect to a linearization close to the pull-back of
$\cO_P(1)$ \cite[\S3]{k2}. Then we have an $SL(2)$-invariant
morphism \[ \phi:\widetilde{P}^s\to \obMz (\pp V,2) \] which
induces an isomorphism
\[ \xymatrix{
\widetilde{Q}:=\widetilde{P}^s/PGL(2)\ar[r]^{\cong} & \obMz (\pp
V,2).}
\]
\end{theorem}

\begin{remark}\lab{rmk1.2}
$P_1^{ss}$ is the locus of points in $P^{ss}$ whose infinitesimal
stabilizer is nontrivial. Hence
$\widetilde{Q}=\widetilde{P}^s/PGL(2)$ is the partial
desingularization defined by Kirwan in \cite{k2}.\end{remark}

\noindent \underline{Proof of Theorem \ref{thm1.1}}. To obtain
$\phi$, we construct a family of stable maps to $\pp V$
parameterized by $\widetilde{P}^s$. The evaluation map
\[
sl(2)^\vee\otimes \Hom (sl(2)^\vee, V)\to V
\]
gives rise to
\[
H^0(\pp^1\times P, \cO(1,1))=sl(2)\otimes \Hom(sl(2)^\vee,V)^\vee
\leftarrow V^\vee=H^0(\pp V,\cO(1)).
\]
Hence we have a rational map
\[
\Phi: \pp^1\times \widetilde{P}^s\to \pp^1\times P \dashrightarrow
\pp V
\]
which is a morphism on $\pp^1\times P_3$ where $P_3=\pp
\mathrm{Hom}_3(sl(2)^\vee, V)$. 
Let us look at the exceptional points of $\Phi$.

Let $\nu$ be a semistable point in $P$. Obviously, when the rank
of $\nu$ is 3, there is no exceptional point. When the rank of
$\nu$ is 2, let $(0,\nu)\in \pp^1\times P^{ss}$ be an
exceptional point. Then we can choose a basis $\{t,s\}$ of
$H^0(\pp^1,\cO(1))$ such that the image of $\nu$ is the span of
$t^2$ and $ts$. But then $\nu$ is a strictly semistable point
whose orbit closure intersects with $P_1$ and thus $\nu$ becomes
unstable after blowing up $P$ along $P_1$ \cite[\S6]{k2}. Hence
such $\nu$ is discarded after the blow-up.

Now suppose the rank of $\nu$ is 1. Then there are two exceptional
points of the rational map $\pp^1\dashrightarrow\pp V$ induced by
$\nu$. But when the two points coincide, $\nu$ becomes unstable
after the blowup along $P_1$ by the same reason as in the previous
paragraph and hence is discarded. Therefore, the exceptional locus
of $\Phi$ consists of two disjoint sections $s_1,s_2$ over
$\widetilde{P}_1^s$, the inverse image of $P_1$ in
$\widetilde{P}^s$.  In particular, the exceptional locus consists
of two codimension 2 subvarieties. We blow up $\pp^1\times
\widetilde{P}^s$ along the exceptional locus and denote the blowup
map by
\begin{equation}\lab{eq2.2}
\xymatrix{ \Gamma\ar[r]^{\rho} & \pp^1\times \widetilde{P}^s
\ar@{-->}[r]^{\Phi} & \pp V}
\end{equation}
and let $D_1,D_2$ be the exceptional divisors. The line bundle
\[
H=\rho^*\cO_{\pp^1\times \widetilde{P}^s}(1,1)\otimes
\cO_{\Gamma}(-D_1-D_2)
\]
is equipped with a surjective morphism $V^\vee\otimes
\cO_\Gamma\to H$ by easy local computation and hence the rational
map \eqref{eq2.2} is a morphism. Therefore, we have a diagram of
morphisms
\[
\xymatrix{ \Gamma\ar[r]^{\Phi\circ \rho}\ar[d]_{\pi_1\circ\rho}
&\pp V\\ \widetilde{P}^s } \] and hence a morphism
$\phi:\widetilde{P}^s\to \obMz(\pp V,2)$. From our construction,
it is evident that $\phi$ is $PGL(2)$-invariant over the dense
open set $P_3$ and thus invariant everywhere.

Since $\pp V$ is convex, $\obMz(\pp V,2)$ is irreducible and
normal. Also by GIT, $\widetilde{Q}$ is projective and hence the
induced map $\widetilde{Q}\to \obMz(\pp V,2)$ is surjective. To
prove that it is an isomorphism it suffices to check injectivity.
It is easy to see by Luna's slice theorem that the exceptional
divisor of $\widetilde{Q}\to Q$ is
a $\pp^{g-2}\times_{\zz_2} \pp^{g-2}$-bundle over $\pp^{g-1}$
which is precisely the set of intersecting lines in $\pp V$. So we
are done. \qed
\bigskip

Note that $\obMz (\pp V,2)=\widetilde{P}^s/PGL(2)=\widetilde{Q}$
is singular along the proper transform of the quotient $Q_2$ of
the locus of rank $\le 2$ homomorphisms $P_2:=\pp
\Hom_2(sl(2)^\vee, V)$ by $PGL(2)$. In fact the proper transform
of $P_2$ in $\widetilde{P}^s$ is the locus of nontrivial
stabilizers which are isomorphic to $\zz_2$. By blowing up along
this singular locus we get a (full) desingularization
$\widehat{Q}$.

\begin{proposition} \lab{prop1.3} $\widehat{Q}$ is the variety of complete
conics $CC(\cB)$ where $\cB\to Gr(3,V)$ is the universal bundle
over the Grassmannian, i.e. $\widehat{Q}$ is the blow-up of $\pp
(S^2\cB)$ along the locus $\pp (S^2_1\cB)$ of rank 1 conics.
\end{proposition}

\begin{proof} For $\alpha\in\Hom (sl(2)^\vee, B)$ for $B\in
Gr(3,V)$, consider the composition
\[ \xymatrix{
B^\vee \ar[r]^{\alpha^\vee} & sl(2)\ar[r] & sl(2)^\vee
\ar[r]^\alpha & B }\] where the middle map is given by the
invariant pairing on $sl(2)$. This gives rise to an isomorphism
\[
\pp \Hom (sl(2)^\vee, \cB)\git PGL(2) \to \pp (S^2\cB)
\]
and hence an isomorphism of blow-ups
\[
bl_{\pp Hom_1(sl(2)^\vee, \cB)}\pp \Hom (sl(2)^\vee, \cB)\git
PGL(2)\cong bl_{ \pp (S^2_1\cB)} \pp (S^2\cB)=CC(\cB).
\]

On the other hand, the inclusion $\cB\hookrightarrow V\otimes
\cO_{Gr(3,V)}$ induces a morphism $\pp \Hom (sl(2)^\vee, \cB)\to
\pp \Hom (sl(2)^\vee, V)$ which preserves the loci of rank $\le i$
homomorphisms for $i=1,2$. Hence we have a map
\[
bl_{\pp Hom_1(sl(2)^\vee, \cB)}\pp \Hom (sl(2)^\vee, \cB)\to
bl_{\pp Hom_1(sl(2)^\vee, V)}\pp \Hom (sl(2)^\vee, V). \] As the
proper transform of $\pp Hom_2(sl(2)^\vee, \cB)$ is a divisor of
the left side, we can further blow up the right side along the
proper transform $\widetilde{\pp \Hom}_2(sl(2)^\vee, V)$ of $\pp
\Hom_2(sl(2)^\vee, V)$ and get a morphism
\[
bl_{\pp Hom_1(sl(2)^\vee, \cB)}\pp \Hom (sl(2)^\vee, \cB)\to
bl_{\widetilde{\pp \Hom}_2(sl(2)^\vee, V)}\left( bl_{\pp
Hom_1(sl(2)^\vee, V)}\pp \Hom (sl(2)^\vee, V) \right) . \] The
quotient of the right side is the blowup $\widehat{Q}$ of
$\widetilde{Q}$ and hence we obtain a map
\[ \phi:CC(\cB)\to \widehat{Q} .\]
We claim this is bijective and therefore an isomorphism. Since
every conic uniquely determines a projective plane $\pp B$
containing it unless it is a double line, $\phi$ is a bijective
over the complement in $\widetilde{Q}$ of the locus of double
lines. Now from the proof of Theorem \ref{thm1.1} the locus of the
double lines in $\widetilde{Q}$ is exactly $\widetilde{\pp
\Hom}_2(sl(2)^\vee, V)\git PGL(2)$ and the normal cones are
$\cc^{g-2}/\zz_2$ by Luna's slice theorem. Hence the fibers of
$\widehat{Q}\to \widetilde{Q}$ over double lines are $\pp^{g-3}$
which parameterizes choices of 3-dimensional subspaces of $V$
containing the image of $sl(2)^\vee$. The fibers of $CC(\cB)\to
\widetilde{Q}$ over a double line corresponds to choices of a
projective plane $\pp B$ containing the line and thus $\pp^{g-3}$.
So $\phi$ is an isomorphism.
\end{proof}

We have blowup maps
\[  \xymatrix{
\pp (S^2\cB) & CC(\cB) \ar[l]_{\Phi_\cB}\ar[r]^{\Phi_{\cB^\vee}} &
\pp (S^2\cB^\vee)}   \] Let $h_1$ be the pull-back to $CC(\cB)$ of
the ample generator $c_1(\cB^\vee)$ of $Pic(Gr(3,V))$; let $h_2$
be the pull-back to $CC(\cB)$ of $c_1(\cO_{\pp (S^2\cB)}(1))$; let
$h_3$ be the class of the exceptional divisor of $\Phi_\cB$. Then
the Picard group of $CC(\cB)$ is generated by $h_1,h_2,h_3$.

The exceptional divisor $\widehat{\Delta}$ of the blowup
$\hQ=CC(\cB)\to \pp(S^2\cB^\vee)$ is $\pp(S^2(\cB^\vee/\cO_{\pp
\cB^\vee}(-1)))$ over $\pp \cB^\vee\cong \pp (S^2_1\cB^\vee)$. Let
$\cA\to Gr(2,V)$ be the universal bundle. Then we have a natural
correspondence
\[\xymatrix{
&\pp \cB^\vee\cong \pp (V/\cA)\ar[dl]\ar[dr]\\
Gr(3,V) && Gr(2,V) }\] As a bundle over $\pp (V/\cA)$, $\pp
(S^2(\cB^\vee/\cO_{\pp \cB^\vee}(-1)))$ is the pull-back of $\pp
(S^2\cA)$ and hence
\begin{equation}\lab{eq3.20}
\widehat{\Delta}\cong \pp (S^2\cA^\vee)\times_{Gr(2,V)}\pp
(V/\cA)\end{equation} Let $\sigma_1$ (resp. $\s_2$) be the
numerical class of a line in a fiber of $\pp (S^2\cA^\vee)$ (resp.
$\pp (V/\cA)$) over $Gr(2,V)$. Let $\s_3$ be the numerical class
of a line in a fiber of $\Phi_\cB$. Then it is straightforward to
compute the pairing of $h_i\cdot \s_j$ for $i,j=1,2,3$ and obtain
the following table:
\begin{equation}\lab{eq3.21}
\begin{matrix}
 & h_1 &h_2&h_3\\
 \s_1&0&1&2\\
 \s_2&1&0&0\\
 \s_3&0&0&-1\end{matrix}\end{equation}
 Also by direct computation the numerical class of
 $\widehat{\Delta}$ is $-2h_1+3h_2-2h_3$. Hence we have
 \begin{equation}\lab{eq3.22}
 \s_1\cdot \widehat{\Delta}=-1,\quad
 \s_2\cdot \widehat{\Delta}=-2,\quad
 \s_3\cdot \widehat{\Delta}=2\end{equation}
As a consequence the restriction of
$\cO_{CC(\cB)}(\widehat{\Delta})$ to a fiber $\pp^2\times
\pp^{g-3}$ of $\widehat{\Delta}\to Gr(2,V)$ is $\cO(-1,-2)$. The
contraction of $CC(\cB)$ along $\pp (S^2\cA^\vee)$ of
$\widehat{\Delta}$ is $\pp(S^2\cB^\vee)$ and
the contraction along $\pp (V/\cA)$ is
$\tQ$. We can further contract $\pp (S^2\cB^\vee)$ along $\pp
(V/\cA)$ and $\tQ$ along $\pp (S^2\cA^\vee)$. They both give
us the same variety, say $\overline{Q}$. In summary we have the
following diagram from contractions in different orders:
\begin{equation}
\xymatrix{ &
\widehat{Q}\ar[dl]_{\pp(V/\cA)}\ar[dr]^{\pp(S^2\cA^\vee)}\\
\tQ\ar[dr]_{\pp(S^2\cA^\vee)} &&\pp (S^2\cB^\vee)
\ar[dl]^{\pp(V/\cA)}\\
&\overline{Q} }\end{equation}

More generally, consider the vector bundle $\cW\to J=Pic^0(X)$ of
rank $g$ given by $\cW=R^1{\pi_J}_*\cL^{-2}(-x)$ where $\pi_J:
J\times X\to J$ is the projection and $\cL\to J\times X$ is
a universal bundle. Let
\[ Q_J:=\pp \Hom(sl(2)^\vee, \cW)\git PGL(2)\]
be the quotient of the projective bundle of homomorphisms of
$sl(2)^\vee$ to a fiber of $\cW\to J$. Let $\tQ_J$ be the blowup of
$Q_J$ along the locus of rank 1 homomorphisms, i.e. the partial
desingularization of $Q_J$. Then analogously as above we obtain
the following.
\begin{theorem}\lab{thm1.2}
$\tQ_J$ is isomorphic to the moduli space $\obMz(\pp\cW/J,2)$ of
relative stable maps of degree 2 to $\pp\cW\to J$.
\end{theorem}
The partial desingularization $\tQ_J$ is singular along the proper
transform $$\Delta_J:=\widetilde{\pp \Hom_2}(sl(2)^\vee, \cW)\git
PGL(2)$$ of the locus of rank $\le 2$ homomorphisms. The blowup
$\hQ_J$ of $\tQ_J$ along $\Delta_J$ is the variety of complete
conics $CC(\cB_J)$ where $\cB_J\to Gr(3,\cW)$ is the universal
rank 3 bundle over the relative Grassmannian. The exceptional
divisor is
\begin{equation}\lab{eq5.131}
\widehat{\Delta}_{J}=\pp(S^2\cA_J^\vee)
\times_{Gr(2,\cW)}\pp(\cW/\cA_J)
\end{equation}
where $\cA_J$ is the universal
rank 2 bundle over the Grassmannian $Gr(2,\cW)$.
Furthermore we have
the following diagram by contractions:
\begin{equation}\lab{eq3.53}
\xymatrix{ &
\widehat{Q}_J=CC(\cB_J)\ar[dl]_{\pp(\cW/\cA_J)}\ar[dr]^{\pp(S^2\cA_J^\vee)}\\
\tQ_J\ar[dr]_{\pp(S^2\cA_J^\vee)} &&\pp (S^2\cB_J^\vee)
\ar[dl]^{\pp(\cW/\cA_J)}\\
&\overline{Q}_J }\end{equation} Theorem
\ref{thm1.2} and \eqref{eq3.53} hold similarly for $\cW_0\to \tX$
instead of $\cW\to J$ where $\cW_0$ is the restriction of $\cW$
via the embedding $\tX\hookrightarrow J$.

A degree $2$ stable map $f:C\to \pp H^1(\xi^{-2}(-x))$ for $\xi\in
J$ gives us a degree $2$ stable map $f:C\to \cN$ by composing with
the embedding $$\pp H^1(\xi^{-2}(-x))\cong \pp Ext^1(\xi,\xi^{-1}(-x))
\hookrightarrow \cN$$ by Proposition \ref{prop2.134}. The
bundle $\pp \cW$ parameterizes extension bundles of a degree 0
line bundle by a degree $-1$ line bundle and thus there is a
natural morphism $\pp \cW\to \cN$. By Theorem \ref{thm1.2} and
this morphism, we obtain the following.
\begin{corollary}\lab{cor1.1} There is an injective morphism
\[ \lambda: \tQ_J\to \obMz(\cN,2) \]
whose image parameterizes degree 2 stable maps of extension type.
\end{corollary}

Injectivity is a consequence of the following result of Narasimhan
and Ramanan.
\begin{proposition} \lab{lem1.5} \cite[6.19]{NR} Let $\xi_1\ne \xi_2\in J$. Then $\pp
H^1(\xi_1^{-2}(-x))$ and $\pp H^1(\xi_2^{-2}(-x))$ intersect if
and only if $(\xi_1,\xi_2)\in \cK_J-\tX$. When they meet, they
intersect at exactly one point transversely.\end{proposition}

In \S6, we will prove that the image of $\lambda$ is normal and
hence $\lambda$ is an isomorphism by Zariski's main theorem.


\section{Hecke correspondence and partial desingularization}

Let $\cM_X=\cM(2,X)$ (resp. $\cM=\cM(2,\cO_X)$, resp. $\cM(2,0)$) be
the moduli space of rank 2 semistable bundles $E$ with $\det
E\cong \cO_X(y-x)$ for some $y\in X$ (resp. $\det E\cong \cO_X$,
resp. $\deg E=0$). Then we have a Cartesian diagram
\[\xymatrix{
\cM(2,\cO_X)\ar@{^(->}[r]\ar[d] & \cM_X\ar[d]\ar@{^(->}[r]&
\cM(2,0)\ar[d]^{\det}\\
\{x\} \ar@{^(->}[r]&X\ar@{^(->}[r]& J}\] where the inclusion
$X\hookrightarrow J$ is $y\mapsto \cO_X(y-x)$.

The projectivization $q_1:\cH\to \cN\times X$ of a universal
bundle over $\cN\times X$ can be thought of as the moduli space of
rank 2 bundles equipped with parabolic structure. Elementary
modification using the parabolic structure gives us a family of
rank 2 semistable bundles on $X$ with determinant $\cO_X(y-x)$ for
some $y\in X$, parameterized by $\cH$ and hence we have a map
$q_0:\cH\to \cM_X$. Therefore we have the Hecke correspondence
\[\xymatrix{ &\cH\ar[dl]_{q_1}\ar[dr]^{q_0}\\
\cN\times X && \cM_X}\]

Although there is no universal bundle for $\cM_X^s$, there is a
projective universal bundle on $\cM_X^s\times X$. The restriction
$\cH_0$ of this $\pp^1$-bundle to $\cM_X^s$ by the graph
\[\alpha:\cM_X^s\to \cM_X^s\times X
\] of the
determinant map $\det:\cM_X^s\to X$ parameterizes a rank 2 stable
bundle $E$ of determinant $\cO_X(y-x)$ and a parabolic structure
at $y$ for $y\in X$. It is elementary to set up moduli problems
for $\cH$ and $\cH_0$ with parabolic bundles and see that $\cH_0$
is the open subscheme $q_0^{-1}(\cM_X^s)$ of $\cH$. Therefore we
have a family of stable maps
\[
\xymatrix{ \cH_0\ar[r]\ar[d] &\cH\ar[r]^{q_1} &\cN\times
X\ar[r]^{\pi_\cN} & \cN\\
\cM_X^s }\] By definition, this is precisely the family of all
Hecke curves in $\cN$ and we obtain a morphism
\[ \psi':\cM_X^s\to \obMz(\cN,2)\]
which is injective by Theorem \ref{thm3.100} (3).
\begin{lemma} $\psi'$ is a birational morphism into an open subvariety
of $\obMz(\cN,2)$.\end{lemma}
\begin{proof} By Theorem \ref{thm3.100} (2), the image of
$\psi'$ is contained in a nonsingular open subvariety of dimension
$3g-2=\dim \cM_X^s$ (Remark \ref{rmk2.4n}). Hence $\psi'$ is an
injective morphism between nonsingular irreducible quasiprojective
varieties of the same dimension. By considering the Jacobian
determinant as in the proof of Theorem
\ref{thm3.18}, $\psi'$ is a birational morphism to an open
subvariety of $\obMz(\cN,2)$.
\end{proof} We claim that $\psi'$ extends to a morphism
$\psi:\tcM_X\to \obMz(\cN,2)$ from the partial desingularization
of $\cM_X$.

Kirwan's partial desingularization has been defined and studied in
\cite{k5}. See also \cite{KL} \S\S 3,5. The singular locus of
$\cM_X$ is the locus of polystable bundles $E=\xi_1\oplus \xi_2$
for $\xi_1,\xi_2\in J=Pic^0(X)$, i.e.
\[ \cK_J=\{ (\xi_1,\xi_2)\in J\times
J\,|\,\xi_1\xi_2=\cO_X(y-x)\text{ for some }y\in X\}/\zz_2 \]
where $\zz_2$ interchanges $\xi_1$ and $\xi_2$. If $\xi_1\ne
\xi_2$, the singularity along $\cK_J$ is by Luna's slice theorem
\[ \cc^{g-1}\oplus \cc^{g-1}\git \cc^*\]
where $\cc^*$ acts with weights $1$ and $-1$ on the summands
respectively. When $\xi_1=\xi_2$, the singularity along $
\tX=\{\xi\in J\,|\, \xi^2=\cO_X(y-x)\text{ for some }y\in X\}$ is
\[ \Hom (sl(2)^\vee,\cW_0)\git PGL(2)\]
where $\cW_0=R^1\pi_*\cL_0^{-2}(-x)$ is a vector bundle of rank
$g$ over $\tX$, $\cL_0$ being the restriction of a Poincar\'e
bundle on $J\times X$ and $\pi:\tX\times X\to \tX$ the projection.
To get the partial desingularization we blow up $\cM_X$ first
along the deepest stratum $\tX$ and then along the proper
transform of the middle stratum $\cK_J$. The exceptional divisor
of the first blowup is
\[\pp\Hom (sl(2)^\vee, \cW_0)\git PGL(2)=Q_{\tX}
\] and its intersection with the proper transform of $\cK_J$ is
the locus of rank 1 homomorphisms
\[\pp\Hom_1 (sl(2)^\vee, \cW_0)\git PGL(2).\] The proper
transform of $Q_{\tX}$ after the second blowup is the partial
desingularization \begin{equation}\lab{eq5.111}
\tQ_{\tX}=bl_{\pp\Hom_1 (sl(2)^\vee, \cW_0)}\pp\Hom (sl(2)^\vee,
\cW_0)\git PGL(2) \end{equation} and the singular locus of
$\tcM_X$ is precisely
\[
\Delta_{\tX}:=bl_{\pp\Hom_1 (sl(2)^\vee, \cW_0)}\pp\Hom_2
(sl(2)^\vee, \cW_0)\git PGL(2)\] and the singularities are
$\zz_2$-quotients. From \S4, $$\Delta_{\tX}\cong
\pp(S^2\cA^\vee_{\tX})\to Gr(2,\cW_0)$$ where $\cA_{\tX}$ is the
universal rank 2 bundle over the Grassmannian  bundle
$Gr(2,\cW_0)$ over $\tX$.
By blowing up $\tcM_X$ once again along $\Delta_{\tX}$ we get a
(full) desingularization $\hcM_X$. As in Proposition
\ref{prop1.3}, the proper transform $\widehat{Q}_{\tX}$ of the
exceptional divisor $Q_{\tX}$ of the first blowup is
\[
CC(\cB_{\tX})\to Gr(3,\cW_0)
\]
and the exceptional divisor of the last blowup $\hcM_{X}\to
\tcM_{X}$ is
\begin{equation}\lab{eq5.130}
\widehat{\Delta}_{\tX}=\pp(S^2\cA^\vee_{\tX})
\times_{Gr(2,\cW_0)}\pp(\cW_0/\cA_{\tX}\oplus \eta)
\end{equation}
for some line bundle $\eta$ on $Gr(2,\cW_0)$.
Let $\widehat{\Sigma}$ be the proper transform of the exceptional
divisor of the second blowup. For $B\in Gr(3,\cW_0)$, let $CC(B)$
denote the fiber of $CC(\cB_{\tX})\to Gr(3,\cW_0)$ over $B$. Then using
the notation of \S4, we have
\[
[\widehat{\Sigma}]|_{CC(B)}=h_3
\]
\[
[\widehat{\Delta}_{\tX}]|_{CC(B)}=-2h_1+3h_2-2h_3
\]
\[
K_{\widehat{Q}_{\tX}}|_{CC(B)}=-(g-4)h_1-6h_2+2h_3
\]
by local computation. (See \cite[\S5]{KL}.) Using \eqref{eq3.21},
we obtain the following table of intersection numbers:
\begin{equation}\lab{eq5.117}
\begin{matrix}
 & \widehat{Q}_{\tX} & \widehat{\Sigma} &\widehat{\Delta}_{\tX} & K_{\hcM_X}\\
 \s_1&0&2&-1&-2\\
 \s_2&1&0&-2&-(g-3)\\
 \s_3&-1&-1&2&-1\end{matrix}\end{equation}
In particular, the intersection numbers of the divisor class
$[\widehat{\Delta}_{\tX}]$ with rational curves $C_1$ in the
fibers of $\pp(S^2\cA^\vee_{\tX})$ and $C_2$ in
$\pp(\cW_0/\cA_{\tX}\oplus \eta)$ are
$$[\widehat{\Delta}_{\tX}]\cdot C_1=-1,
\qquad [\widehat{\Delta}_{\tX}]\cdot C_2=-2.$$ Hence the
restriction of the normal bundle
$\cO_{\hcM_X}(\widehat{\Delta}_{\tX})$ to the fibers
$\pp^2\times \pp^{g-2}$ of $\widehat{\Delta}_{\tX}\to Gr(2,\cW_0)$
is $\cO(-1,-2)$. The contraction of $\hcM_X$ along
$\pp(\cW_0/\cA\oplus \eta)$ gives us the orbifold $\tcM_X$ while
the contraction of $\pp(S^2\cA^\vee_{\tX})$ is a nonsingular
variety $\cM_X^\#$. We will prove in \S7 that this nonsingular
variety is one of the two irreducible components of the Hilbert
scheme $\bH$ of conics in $\cN$. We can further contract
$\pp(\cW_0/\cA_{\tX}\oplus \eta)$ to obtain an orbifold
$\cM_X^\flat$ which turns out to be an irreducible component of
the Chow scheme of conics in $\cN$. In summary we have the
following commutative diagram whose maps are all contractions (in
different orders).
\begin{equation}
\xymatrix{ &
\widehat{\cM}_X\ar[dl]_{\pp(\cW_0/\cA_{\tX}\oplus \eta)}
\ar[dr]^{\pp(S^2\cA^\vee_{\tX})}\\
\tcM_X\ar[dr]_{\pp(S^2\cA^\vee_{\tX})}
&&\cM_X^\# \ar[dl]^{\pp(\cW_0/\cA_{\tX}\oplus \eta)}\\
&\cM_X^\flat }\end{equation}

In \cite{CCK} \S\S5,6, we constructed a family of Hecke cycles
parameterized by (a neighborhood of each point in) the
desingularization $\hcM(2,\cO_X)$ of $\cM(2,\cO_X)$ and during the
proof a family of maps parameterized by the partial
desingularization $\tcM(2,\cO_X)$. It is straightforward to modify
the constructions to obtain the following.
\begin{theorem}\lab{thm5.90} The birational morphism $\psi':\cM_X^s\to
\obMz(\cN,2)$ extends to an injective morphism
\begin{equation}\lab{eq5.90}
\psi:\tcM_X\to \obMz(\cN,2).\end{equation} The image of $\psi$ is
an irreducible component parameterizing all stable maps of the
following types: \begin{enumerate} \item[(i)] Hecke curves
$f_{E,y}$ for $E\in \cM_X^s$, $\det E\cong \cO_X(y-x)$ \item[(ii)]
degree 2 maps $f:C\to \pp H^1(\xi^{-2}(-x))$ for $\xi\in \tX$
\item[(iii)] unions of two intersecting lines $l_1\cup l_2$ for
$l_1\subset \pp H^1(\xi_1^{-2}(-x))$, $l_2\subset \pp
H^1(\xi_2^{-2}(-x))$, $(\xi_1, \xi_2)\in
\cK_J-\tX$.\end{enumerate}\end{theorem}
\begin{proof} We leave the modification of \cite{CCK} \S\S5,6 to
the reader. The local computation in the proofs of Lemma 5.2 and
Lemma 6.1 shows that the stable maps parameterized by $\tcM_X$ are
of the above types. \end{proof}

In \S6, we will prove that the image of $\psi$ is normal and thus
$\psi$ is an isomorphism onto an irreducible component of
$\obMz(\cN,2)$.


\section{Moduli space of stable maps of degree 2 to $\cN$}

In this section we prove Theorem \ref{thm0.1}. By Proposition
\ref{prop4.1}, there are only 3 types of stable maps of degree 2.
\begin{lemma}\lab{lem5.0} A stable map of degree 2 to $\cN$ is
one of the following:
\begin{enumerate}
\item[(i)] Hecke curves $f_{E,y}$ for $E\in \cM_X^s$, $\det E\cong
\cO_X(y-x)$ \item[(ii)] degree 2 maps $f:C\to \pp
H^1(\xi^{-2}(-x))$ for $\xi\in J=Pic^0(X)$ \item[(iii)] union of
two intersecting lines $l_1\cup l_2$ for $l_1\subset \pp
H^1(\xi_1^{-2}(-x))$, $l_2\subset \pp H^1(\xi_2^{-2}(-x))$,
$\xi_1\ne \xi_2\in J$.\end{enumerate}
\end{lemma}

By Proposition \ref{lem1.5}, if $\pp H^1(\xi_1^{-2}(-x))$ and $\pp
H^1(\xi_2^{-2}(-x))$ intersect, then $(\xi_1,\xi_2)\in \cK_J-\tX$.
Hence stable maps of type (iii) are all contained in the image of
$\psi:\tcM_X\to \obMz(\cN,2)$ by Theorem \ref{thm5.90}.

We proved in \S5 that there is an injective morphism
$\psi:\tcM_X\to \obMz(\cN,2)$ whose image is an irreducible
component $\bM_H$, which we call the \emph{Hecke component},
parameterizing all stable maps of types (i), (iii) and stable maps
of type (ii) with $\xi\in \tX$. We also proved in \S4 that there
is an injective morphism $\lambda:\tQ_J\to \obMz(\cN,2)$ whose
image parameterizes stable maps of type (ii) with $\xi$ not
necessarily in $\tX$. Since $\tQ_J$ is irreducible projective and
$\lambda$ is an injective map whose image contains the complement
of the Hecke component, $\lambda$ is bijective onto an irreducible
component $\bM_E$, which we call the
\emph{extension component}. Obviously the intersection of the two
components parameterizes degree 2 maps $f:C\to \pp
H^1(\xi^{-2}(-x))$ for $\xi\in \tX$. Hence the intersection is the
image of $\tQ_{\tX}$ by both $\lambda$ and $\psi$.

The purpose of this section is to prove the following.
\begin{proposition}\lab{prop5.100} The Hecke component $\bM_H$
and the extension component $\bM_E$  are normal varieties
intersecting transversely along a subvariety isomorphic to
$\tQ_{\tX}$.\end{proposition}

\noindent \underline{Proof of Theorem \ref{thm0.1}.} As $\tcM_X$
and $\tQ_J$ are irreducible normal projective, $\lambda$ and
$\psi$ are bijective morphisms onto normal varieties and thus
birational (see the proof of Theorem \ref{thm3.18}). By Zariski's
main theorem, $\psi$ and $\lambda$ are isomorphisms. \qed

\bigskip
The rest of this section is devoted to a proof of Proposition
\ref{prop5.100}. We already know that $\obMz(\cN,2)$ is normal at
Hecke curves (Remark \ref{rmk2.4n}).

To prove normality at the stable maps of type (ii) and (iii) of
Lemma \ref{lem5.0}, we need to know the restriction of the tangent
bundle $T\cN$ of $\cN$ to rational curves in $\cN$ of extension
type.

\begin{proposition} \lab{prop5.10}
(1) Let $l$ be a line in $\pp H^1(\xi^{-2}(-x))$ for $\xi\in J$.
Then we have
\begin{equation}\lab{eq5.100}
T\cN|_l\cong \left\{\begin{array}{lll} \cO(2)\oplus
\cO(1)^{g-2}\oplus \cO(-1)^{g-2}\oplus \cO^g & \text{if}&
\xi\notin
\tX\\
\cO(2)\oplus \cO(1)^{g-1}\oplus \cO(-1)^{g-1}\oplus \cO^{g-2}
&\text{if}& \xi\in \tX.
\end{array}\right.
\end{equation}

(2) If $f:\pp^1\to \pp H^1(\xi^{-2}(-x))$ for $\xi\in J$ is of
degree 2, then
\begin{equation}\lab{eq5.101}
f^*T\cN\cong \left\{\begin{array}{lll} \cO(4)\oplus
\cO(2)^{g-2}\oplus \cO(-2)^{g-2}\oplus \cO^g & \text{if}&
\xi\notin
\tX\\
\cO(4)\oplus \cO(2)^{g-1}\oplus \cO(-2)^{g-1}\oplus \cO^{g-2}
&\text{if}& \xi\in \tX.
\end{array}\right.
\end{equation}
\end{proposition}
\begin{proof}
For a line $l\subset \pp H^1(\xi^{-2}(-x))\cong \pp^{g-1}$,
$\xi\in J$, the tangent bundle of $\pp H^1(\xi^{-2}(-x))$
restricted to $l$ is
\begin{equation}\lab{eq5.13} \cO(2)\oplus
\cO(1)^{g-2}\end{equation} The computation of the normal bundle of
$\pp H^1(\xi^{-2}(-x))$ in $\cN$ is similar to the proof of Lemma
\ref{lem2.117}. We consider the universal extension bundle $\cU$
over $\pp H^1(\xi^{-2}(-x))\times X$:
\[ 0\to\pi_X^*\xi^{-1}(-x)\otimes \cO_{\pp}(1)\to \cU\to
\pi_X^*\xi\to 0\] By \cite[6.7]{NR}, the tangent bundle $T\pp$ of
$\pp H^1(\xi^{-2}(-x))$ fits into an exact sequence
\begin{equation}\lab{eq5.10}
0\to T\pp\to R^1{\pi_\pp}_*(\cU\otimes \pi_X^*\xi^{-1})\to
R^1{\pi_\pp}_*\cO=H^1(X,\cO_X)\otimes \cO_{\pp}\to 0\end{equation}
where $\pi_\pp$ denotes the projection onto $\pp H^1(\xi^{-2}(-x))$.
We also have an exact sequence
\begin{equation}\lab{eq5.11} 0\to \cU\otimes
\pi_X^*\xi^{-1}\to End_0\cU\to \pi_X^*\xi^2(x)\otimes \cO_\pp
(-1)\to 0.
\end{equation}
This gives us a diagram
\begin{equation}\lab{eq5.12}\xymatrix{ 0\ar[r]
&H^0(\xi^2(x))\otimes \cO_\pp(-1)\ar@{=}[d] \ar[r]& R^1{\pi_\pp}
_*(\cU\otimes \pi_X^*\xi^{-1})\ar[r]\ar[d] & R^1{\pi_\pp}
_*End_0\cU\ar@{->>}[d]\\
&H^0(\xi^2(x))\otimes \cO_\pp(-1)\ar[r] & H^1(\cO_X)\otimes \cO_\pp &
H^1(\xi^2(x))\otimes \cO_\pp(-1).}
\end{equation}
 Suppose
$\xi\notin\tX$. Then $H^0(\xi^2(x))=0$ and
$H^1(\xi^2(x))=\cc^{g-2}$. From \eqref{eq5.13}, \eqref{eq5.10} and
\eqref{eq5.12}, we deduce that
\[
T\cN|_l\cong \cO(2)\oplus \cO(1)^{g-2}\oplus \cO(-1)^{g-2}\oplus
\cO^g
\] when $\xi\notin \tX$.

Now suppose $\xi\in \tX$. Then $H^0(\xi^2(x))=\cc$ and
$H^1(\xi^2(x))\cong \cc^{g-1}$. Since an extension bundle is
uniquely determined by the epimorphism $E\to \xi$, we have
\[ {\pi_{\pp}}_*Hom (\pi_X^*\xi,\pi_X^*\xi)\cong
{\pi_\pp}_*Hom(\cU, \pi_X^*\xi)\] and thus the bottom horizontal
in \eqref{eq5.12} is a monomorphism. Because
\[  H^1(\cO)\otimes \cO_{\pp}/H^0(\xi^2(x))\otimes \cO_\pp (-1) \cong
\cO(1)\oplus \cO^{g-2}\] we deduce that
\[
T\cN|_l\cong \cO(2)\oplus \cO(1)^{g-1}\oplus \cO(-1)^{g-1}\oplus
\cO^{g-2}
\] when $\xi\in \tX$.

By the same computation for conics, we obtain \eqref{eq5.101}.
\end{proof}

We prove normality of $\obMz(\cN,2)$ at the stable maps of type
(iii) in Lemma \ref{lem5.0}. Let $f_1,f_2:\pp^1\to \cN$ be the
lines $l_1, l_2$ with the intersection point $\theta$. Let
$f:\pp^1\cup \pp^1\to \cN$ be the morphism from the gluing of the
points corresponding to $\theta$. From the exact sequence,
\[
0\to f^*T\cN\to f_1^*T\cN\oplus f_2^*T\cN\to T\cN|_\theta \to 0
\]
we get the exact sequence
\[
\xymatrix{ 0\ar[r] & H^0(f^*T\cN) \ar[r] & H^0(f_1^*T\cN)\oplus
H^0(f_2^*T\cN) \ar[r] & H^0(T\cN|_\theta)\ar[r]&}\]
\[\xymatrix{{}\ar[r]&H^1(f^*T\cN)\ar[r] &H^1(f_1^*T\cN)\oplus
H^1(f_2^*T\cN)\ar[r] & 0}
\]
By Proposition \ref{prop5.10} (1), $H^1 (f_i^*T\cN)=0$ and
$H^0(f^*_iT\cN)\cong \cc^{3g-1}$. Hence normality at $[f]$ is a
consequence of the following.
\begin{lemma}
$H^0(f_1^*T\cN)\oplus H^0(f_2^*T\cN) \to H^0(T\cN|_\theta)$ is
surjective.
\end{lemma}
\begin{proof} Because $\xi_1\ne \xi_2$ and
$\xi_1\xi_2=\cO_X(y-x)$, we have $\xi_1,\xi_2\notin \tX$. From
\eqref{eq5.10} and \eqref{eq5.12}, we have an isomorphism
\begin{equation}\lab{eq5.20}
T\cN|_\theta =T\pp H^1(\xi_1^{-2}(-x))|_\theta \oplus
H^1(\cO)\oplus H^1(\xi_1^2(x))
\end{equation}
The three summands are respectively the subspaces of positive,
zero and negative degrees for $T\cN|_{l_1}$. Note that $T\pp
H^1(\xi_1^{-2}(-x))|_\theta=H^1(\xi_1^{-2}(-x))/\cc e_1$ where
$e_1$ is the extension class of
\[
0\to \xi_1^{-1}(-x)\to E\to \xi_1\to 0
\]
at $\theta=[E]$. We also have an exact sequence
\[
0\to \xi_2^{-1}(-x)\to E\to \xi_2\to 0
\]
and $\xi_1(-y)\cong \xi_2^{-1}(-x)$. Thus $\cc e_1$ is the kernel
of
\[
H^1(\xi_1^{-2}(-x))=Ext^1(\xi_1,\xi_1^{-1}(-x))\to
Ext^1(\xi_2^{-1}(-x),\xi_1^{-1}(-x))\cong H^1(\xi_1^{-1}\xi_2).
\]
Therefore, $T\pp H^1(\xi_1^{-2}(-x))|_\theta\cong
H^1(\xi_1^{-1}\xi_2)=H^1(\xi_2^2(x-y)).$ From the exact sequence
\[
0\to \xi_2^2(x-y)\to \xi_2^2(x)\to \cO_y\to 0
\]
we get an exact sequence
\[
0\to \cc\to H^1(\xi_2^2(x-y))\to H^1(\xi_2^2(x))\to 0
\]
Hence \eqref{eq5.20} is now
\[
T\cN |_\theta =H^1(\xi_2^{2}(x))\oplus \cc \oplus H^1(\cO)\oplus
H^1(\xi_1^2(x))
\]
The first three summands are the subspaces of nonnegative degrees
for $T\cN |_{l_1}$ while the last three summands are the subspaces
of nonnegative degrees for $T\cN |_{l_2}$. This implies the desired
surjectivity. \end{proof}

Now we consider the stable maps $f:C\to \cN$ of type (ii) in Lemma
\ref{lem5.0}. Suppose $\xi\notin \tX$. These stable maps are
parameterized by the complement
\[
\obMz(\cN,2)-\psi (\tcM_X)
\]
which is an open set since $\tcM_X$ is projective. By
\eqref{eq5.101}, $$h^1(f^*T\cN )=5+3(g-2)+g=4g-1$$ and hence the
dimension of the deformation space at $[f]\in \obMz(\cN,2)$ is
constantly $4g-4$ which is exactly the dimension of $\tQ_J$. Hence
the dimension of $\obMz(\cN,2)-\psi(\tcM_X)$ is at most $4g-4$. But
we have a bijective morphism
\begin{equation}\lab{eq5.110}
\lambda:\tQ_J-\tQ_{\tX}\to \obMz(\cN,2)-\psi(\tcM_X)
\end{equation}
By considering the Jacobian determinant, we deduce that
$\obMz(\cN,2)-\psi(\tcM_X)$ is irreducible of dimension exactly
$4g-4$ and hence it has at worst quotient singularities. In
particular, $\obMz(\cN,2)-\psi(\tcM_X)$ is normal.

Finally let $f:C\to \cN$ be a stable map of type (ii) with $\xi
\in \tX$.  Then $f$ is in the image of $\tQ_{\tX}$ both by $\psi$
and by $\lambda$ if $\tQ_{\tX}$ is considered as a divisor of
$\tcM_X$ by \eqref{eq5.111} and as a subvariety of $\tQ_J$ by the
embedding $\tX\hookrightarrow J$. Hence $f$ lies in the
intersection of the two components $\bM_H$ and $\bM_E$. By
\eqref{eq5.101}, $h^0(f^*T\cN)=4g$. The proof of Proposition
\ref{prop5.100} is now complete if the deformation of $f$ in
$\bM_H$ is cut out by $(g-1)$ linear equations and the deformation
of $f$ in $\bM_E$ is cut out by $1$ linear equation, transversal
to the $(g-1)$-equations, on $H^0(f^*T\cN)$.

Let us consider the deformation of $f$ in $\bM_E$. By assumption
$f$ factors through the projective bundle $\pp \cW$ over $J$. The
deformation of $f$ as a map to $\pp \cW$ is from \eqref{eq5.11}
$$H^0(R^1{\pi_{\pp}}_*(\cU\otimes \pi_X^*\xi^{-1}))$$ using the
notation of the proof of Proposition \ref{prop5.10}. By taking
cohomology of \eqref{eq5.12} we obtain a diagram
\small
\begin{equation}\lab{eq5.121}
\xymatrix{  0\ar[d] &0\ar[d]\\
H^0(f^*T\pp H^1(\xi^{-2}(-x)))\ar[d]\ar@{=}[r] &H^0(f^*T\pp
H^1(\xi^{-2}(-x)))\ar[d] \\
H^0(R^1{\pi_{\pp}}_*(\cU\otimes \pi_X^*\xi^{-1}))
\ar@{^(->}[r]\ar[d] & H^0(f^*T\cN)\ar[d]\ar@{->>}[r]
&H^1(\cO_{\pp^1}(-2))=\cc\ar@{=}[d] \\
H^1(\cO_X)\ar@{^(->}[r]\ar[d] & H^0([H^1(\cO_X)\otimes
\cO_{\pp^1}]/\cO_{\pp^1}(-2))\ar@{->>}[r]\ar[d] & H^1(\cO_{\pp^1}(-2))=\cc\\
 0& 0 }
\end{equation}
\normalsize Therefore the deformation in $\bM_E$ is cut out by one
linear equation given by the surjective linear map
$H^0(f^*T\cN)\to H^1(\cO_{\pp^1}(-2))=\cc$ in \eqref{eq5.121}.
This implies the normality of the extension component $\bM_E$.

Let us consider the deformation of $f$ in $\bM_H$. Let $\tau$ be
$\psi^{-1}(f)\in \tcM_X$. The embedding $\tX\hookrightarrow J$
gives us a subspace $\cc\hookrightarrow H^1(\cO_X)$. The
deformation space $V_Q$ of $f$ in the divisor $\tQ_{\tX}$ of
$\tcM_X$ fits into a diagram \small
\begin{equation}\lab{eq5.122}
\xymatrix{0 \ar[d]&  0\ar[d] &0\ar[d]\\
H^0(f^*T\pp H^1(\xi^{-2}(-x)))\ar@{=}[r]\ar[d] &H^0(f^*T\pp
H^1(\xi^{-2}(-x)))\ar[d]\ar@{=}[r] &H^0(f^*T\pp
H^1(\xi^{-2}(-x)))\ar[d] \\
V_Q\ar[r]\ar[d] &H^0(R^1{\pi_{\pp}}_*(\cU\otimes \pi_X^*\xi^{-1}))
\ar@{^(->}[r]\ar[d] & H^0(f^*T\cN)\ar[d]\\
\cc\ar[r]\ar[d] &H^1(\cO_X)\ar@{^(->}[r]\ar[d] &
H^0([H^1(\cO_X)\otimes
\cO_{\pp^1}]/\cO_{\pp^1}(-2))\ar[d] \\
 0& 0 &0}
\end{equation}
\normalsize The dimension of $V_Q$ is at most
$h^0(\cO_{\pp^1}(4)\oplus \cO_{\pp^1}(2)^{g-2})+1=3g$. To find the
image of the deformation of $\tau$ in the normal direction of the
divisor $\tQ_{\tX}$ in $\tcM_X$ we look at the local computation
in \cite[\S6]{CCK}. By \cite{CCK} (6.5), if we deform $f$ (or
$\tau$) in the normal direction of $\tQ_{\tX}$ in $\tcM_X$, the
diagonal part of the transition matrices gives us a normal vector
field ($z$ in the notation of \cite{CCK}) of the divisor
$\tQ_{\tX}$ times the tangent vector field ($-a_{ij}-tc_{ij}$ in
\cite{CCK}) of the conic ($b_{ij}-2ta_{ij}-t^2c_{ij}$ in
\cite{CCK}). In particular, it gives us $\cO_{\pp^1}\subset
[H^1(\cO_X)\otimes \cO_{\pp^1}]/\cO_{\pp^1}(-2)$ since the
pull-back of the normal bundle of $\tQ_{\tX}$ by $f$ is
$\cO_{\pp^1}(-2)$ by \eqref{eq5.117}. As the normal deformation is
given by a nonconstant vector field ($-a_{ij}-tc_{ij}$ in
\cite{CCK}), its image of the linear equation $H^0(f^*T\cN)\to\cc$
for $\bM_E$ is nonzero (see \eqref{eq5.121}). In summary, the
deformation of $f$ in $\bM_H$ is given by an extension of $\cc$ by
$V_Q$ which is transversal to $\bM_E$ and is of codimension at
least $g-1$. The proof is similar when $C$ is reducible.
Therefore, $\bM_H$ is normal everywhere and we finished a proof of
Proposition \ref{prop5.100}.

\begin{remark}
If one wants to find the moduli theoretic meaning of the partial
desingularization of $\cM=\cM(2,\cO_X)$ instead of $\cM_X$, we
just have to consider relative stable maps of degree $2$
$$\xymatrix{\Gamma\ar[rr]^f\ar[dr] &&\cN_X\ar[dl]^{\det}\\
&X}$$ where $\Gamma$ is a flat family of semistable curves of
genus $0$ over $X$ and $\cN_X$ is the moduli space of stable
bundles of rank 2 with determinant $\cO_X(-y)$ for some $y\in X$.
We leave the details to the reader.
\end{remark}

\section{Hilbert scheme and Chow scheme}

In this section we describe the Hilbert scheme $\bH$ and the Chow
scheme $\bC$ of \emph{conics} in $\cN$, i.e. subschemes in $\cN$
with Hilbert polynomial $2m+1$. Let $\bM=\obMz(\cN,2)$ and
$\widehat{\bM}$ be the blowup of the disjoint union of the two
components of $\bM$ along the singular locus ${\Delta}_J\cup
{\Delta}_{\tX}$.

We recall the following.
\begin{lemma} \cite[4.2,4.3]{NR} Let $C\subset \cN$ be a conic.
Then $C$ is isomorphic to a conic in $\pp^2$ or a total thickening
of a line $\pp^1$ in $\pp^2$. \end{lemma}

Therefore $\bH$, $\bC$ and $\bM$ are all isomorphic over the
complement of the locus ${\Delta}_J\cup {\Delta}_{\tX}$ of double
covers of a line in $\cN$. Thus we have a diagram
\begin{equation}
\xymatrix{ & \widehat{\bM}\ar[dl]
\ar@{-->}[dr]\\
\bM\ar[dr]
&&\bH\ar[dl]\\
&\bC }\end{equation} In \cite[\S\S5,6]{CCK}, we constructed a
family of Hecke cycles parameterized by (a neighborhood of each
point of) the (full) desingularization of $\hcM$ and thus obtained
a morphism from $\hcM$ to the Hilbert scheme of Hecke cycles.
Further we proved the map is a blowup of the Hilbert scheme of
Hecke cycles. It is straightforward to modify the constructions of
\cite{CCK} to prove the following.
\begin{theorem}\lab{thm7.1} There is a morphism
$\widehat{\bM}\to \bH$ extending the birational map
$\widehat{\bM}\dashrightarrow \bH$. Furthermore, the morphism is
the contraction of $\pp^2(S^2\cA_J^\vee)$ via \eqref{eq5.131},
\eqref{eq5.130}.
\end{theorem}
Next let us consider the natural map from the Hilbert scheme $\bH$
to the Chow scheme $\bC$ (\cite{Kollar}). This is not an
isomorphism along the image of $\widehat{\Delta}_J\cup
\widehat{\Delta}_{\tX}$. The points in $Gr(2,\cW)$ determine a
line $l$ in $\pp \cW$ and the fibers of $\pp (\cW/\cA_J)\to
Gr(2,\cW)$ parameterize a thickening of $l$ by $\cO(-1)$. Hence
the morphism $\bH\to \bC$ factors through the contraction of $\pp
(\cW/\cA_J)$ for the extension component and of $\pp
(\cW/\cA_{\tX}\oplus \eta)$ for the Hecke component. The resulting
map is one-to-one and the Chow scheme is seminormal. Therefore we
deduce that the contraction of $\bH$ is isomorphic to the Chow
scheme $\bC$. So we obtain the following.
\begin{theorem}\lab{thm7.2} (1) The Hilbert scheme $\bH$ of conics
in $\cN$ consists of two irreducible nonsingular components,
intersecting transversely.

(2) $\bM, \bH, \bC$ are related by \eqref{eq0.6} which is the
diagram of contractions of $\widehat{\Delta}_J\cup
\widehat{\Delta}_{\tX}$ in different orders.
\end{theorem}

\bibliographystyle{amsplain}

\end{document}